\RequirePackage{fix-cm}
\documentclass[smallextended,draft]{svjour3}       % onecolumn (second format)
\smartqed  % flush right qed marks, e.g. at end of proof
\usepackage{graphicx}
 \usepackage{mathptmx}      % use Times fonts if available on your TeX system
\usepackage{amsfonts,latexsym,amsmath,amssymb,mathrsfs}
\usepackage{caption}
\usepackage{yfonts}
\usepackage[dvipsnames]{xcolor}
\usepackage{cite}
\usepackage{enumitem}

\newtheorem{assumption}{Assumption}

\DeclareMathOperator{\diff}{d\!} 
 \journalname{Mathematics of Control, Signals, and Systems}
\begin{document}

\title{Coercive ISS-Lyapunov functionals for regular infinite-dimensional systems and applications\thanks{This work has been supported by the Agence National de la Recherche (ANR) under the ROTATION project (grant no. ANR-24-CE48-0759).}
}
%\subtitle{Do you have a subtitle?\\ If so, write it here}

%\titlerunning{Short form of title}        % if too long for running head

\author{Swann Marx}

%\authorrunning{Short form of author list} % if too long for running head

\institute{
S. Marx \at
        LS2N, \'Ecole Centrale de Nantes $\&$ CNRS UMR 6004, F-44000 Nantes, France.    \\
     %         Tel.: +123-45-678910\\
      %        Fax: +123-45-678910\\
              \email{swann.marx@ls2n.fr}           %  \\
%             \emph{Present address:} of F. Author  %  
}

\date{Received: \today / Accepted: date}
% The correct dates will be entered by the editor

\maketitle

\begin{abstract}
This paper proposes the construction of a coercive ISS-Lyapunov functional for linear regular infinite-dimensional system. Indeed, as already known, Lyapunov functionals for infinite-dimensional systems might be not coercive. Under the assumption that there exists an exactly observable output, we are able to make coercive a Lyapunov functional which is not coercive under additional regularity assumption. We discuss also about the potential applications of such a Lyapunov functional in singular perturbation theory and output regulation. The results are illustrated on a non-trivial equation, namely, the Korteweg-de Vries equation.
\end{abstract}

\section{Introduction}

This paper is focused on the construction of ISS Lyapunov functional for regular linear infinite-dimensional systems. Regular linear infinite-dimensional systems are well-posed linear systems (in the sense given by \cite{staffans2005well,weiss1994regular}) for which an additional regularity is supposed for its related transfer function. It turns out that many infinite-dimensional systems satisfy such a property, see e.g., \cite{tucsnak2014well} for a list of such systems. We may for instance mention parabolic equations that satisfy such a property, and we will see that the Korteweg-de Vries equation (studied in \cite{cerpa2013control}) satisfies this property. This additional property allows us to extend many well-known results in finite-dimensional, as is well explained in \cite{weiss1994regular}.

In parallel, ISS theory has been developed for infinite-dimensional systems \cite{Mironchenko2020Input,prieur2011iss,jacob2020noncoercive,Mironchenko2018Characterizations,mironchenko2023lyapunov}, while it is nowadays well-known in the finite-dimensional context \cite{iss_sontag,jiang1994small,Jiang2001Input} and used in many various contexts such as observer design or output regulation. As is well explained in \cite{jacob2020noncoercive}, the ISS property is closely related to the well-posedness of the infinite-dimensional system under consideration: indeed, in \cite[Theorem 5.3.]{jacob2020noncoercive}, it is explained that as soon as a Lyapunov functional exists for the system without input, then the admissibility of the control operator is sufficient to ensure the ISS of the linear infinite-dimensional system. The admissibility of control operators is crucial to establish well-posedness of infinite-dimensional linear systems admitting inputs \cite{weiss1989admissibility,tucsnak2009observation}. In a sense, we are following the same path in this paper, but we add in the description of our system an output, making the analysis of the latter more difficult.

Another point to be discussed here is the lack of coercivity of Lyapunov functionals for linear infinite-dimensional systems that are exponentially stable. This result has been established, for instance, in \cite{hante2011converse,Mironchenko2018Characterizations,jacob2020noncoercive}. The lack of coercivity means that, even if a linear system is exponentially stable, it might happen that the Lyapunov functional under consideration is not equivalent to the usual norm of the state space. This lack of coercivity might be an issue when considering output regulation problems as the ones presented in \cite{terrand2019adding,paunonen2010internal}. 

In \cite[Remark 4.2.]{lorenzetti2023saturating}, a very nice technique has been used in order to render a non-coercive Lyapunov functional coercive. It works as follows: if a linear system admits an infinite-time admissible and exactly observable output, then one can make a non-coercive Lyapunov functional coercive by adding the observability Gramian in the definition of the Lyapunov functional. Because of the properties of the output, the observability Gramian is equivalent to the norm in the state space. 

This article relies on this trick, but furthermore adds input in the definition of the system. Therefore, from a Lyapunov functional that is a non-coercive ISS Lyapunov functional, we can make the latter coercive by adding another term (related to Gramian observability operator) under the assumption that the triple $(A,B,C)$ is regular in the sense given by \cite{curtain1989well,weiss1994transfer,weiss1997dynamic}. Moreover, in contrast with \cite{lorenzetti2023saturating}, we keep the influence of the output in the Lyapunov functional. Finally our Lyapunov functional $V$ satisfies the following inequality (if we consider strong solutions in a sense that will be explained later on):
\begin{equation}
    \frac{\diff V}{\diff t}(z) \leq -\alpha_1 \Vert z\Vert^2_H + \alpha_2 \Vert u\Vert^2_U - \alpha_3 \Vert Cz\Vert^2_H,
\end{equation}
 with $\alpha_1,\alpha_2,\alpha_3>0$ and where $u\in U$ is an input and $z$ is the state of the system. The fact that the output appears as in the latter inequality looks like the strict output passivity discussed in \cite{marx2025impedance}. Such a Lyapunov functional has been used widely in \cite{balogoun2021iss,marx2024singular} for the linearized version of the Korteweg-de Vries equation. This passivity property will be very useful for stability purposes. 

Indeed, in the second part of this paper, we will discuss the potential applications of such a property. First, we will show that, under some structural assumptions, a singular perturbation theorem can be stated. Roughly speaking, the singular perturbation technique (see e.g., \cite{kokotovic1999singular}) concerns coupled systems admitting different time-scales, i.e. a fast system and slower one. The singular perturbation method consists in decoupling the coupled system into two systems: the first one (called the reduced order system) represents the slow dynamics, and the second (called the boundary layer system) describes the fast dynamics. The singular perturbation method states that, as soon as both subsystems are stable, then the full system is also stable provided that the fast system is sufficiently fast. This method has been extended to the infinite-dimensional case \cite{tang2017stability,cerpa-prieur2017,marx2024singular,arias2023frequency,arias2025stability}, but it turns out that some counter-examples have been found in \cite{tang2017stability,cerpa-prieur2017}. We will provide a positive answer to this question in the case of a fast infinite-dimensional systems is coupled with a slow ordinary differential equation. From this result, we will show that, in the spirit of \cite{lorenzetti2022saturating}, we can also achieve output regulation of linear systems using the result dealing with singular perturbation.

The paper is divided as follows: in Section \ref{sec:coercive}, we discuss about regular linear system and Lyapunov functionals, and we state and prove our main result about making coercive an ISS non-coercive Lyapunov functional. In Section \ref{sec:applications}, we talk about the potential applications existing when assuming the existence of such a coercive Lyapunov functional: first, we discuss the singular perturbation method (i.e., coupled systems with different time-scales), and second, we talk about the design of a PI controller for an infinite-dimensional system. Section \ref{sec:conclusion} proposes some concluding remarks and discusses about further research lines to follow.  
\section{Construction of a coercive Lyapunov functional}

\label{sec:coercive}

\subsection{Regular linear systems}

We consider $H,U,Y$ three real Hilbert spaces, where $H$ is the state space, $U$ the input state space and $Y$ the output space. In the sequel, given any Hilbert space $W$, we denote by $\mathrm I_W$ the identity operator of this space, and by $\mathcal L(W_1,W_2)$ the set of bounded (or continuous) operators from $W_1$ to $W_2$, with $W_1,W_2$ two real Hilbert spaces. We also have a triple of operators $(A,B,C)$ which defines the following dynamical system:
\begin{equation}
\label{eq:ISS}
\left\{
\begin{aligned}
&\frac{\diff}{\diff t} z(t) = Az(t) + Bu(t),\\
& y(t) = C z(t),\\
& z(0) = z_0.
\end{aligned}
\right.
\end{equation}
where $A:D(A)\subset H \rightarrow H$, with $D(A)$ densely defined in $H$, and A generates a strongly continuous semigroup $(\mathbb T_t)_{t\geq 0}$, $B\in \mathcal{L}(\mathbb R,H_{-1})$, and $C\in \mathcal{L}(H_1,\mathbb R)$, where $H_{-1}$ is defined as the completion of $H$ with respect to the norm $\Vert (s\mathrm I_H-A)^{-1}z\Vert_H$ with $z\in H$ and $s\in \rho(A)$, and $H_1$ is the space $D(A)$ equipped with the norm $\Vert (s\mathrm I_H-A)z\Vert_H$ with $z\in D(A)$ and $s\in \rho(A)$. The input $u$ is such that $u\in L_{\mathrm loc}^2([0,\infty);U)$. The construction of such a Gelfand triple $(H_{-1},H,H_1)$ is nowadays quite standard and can be found, for example, in \cite{tucsnak2009observation}. We suppose that $B$ and $C$ are admissible operators for $\mathbb T$. Since $C\in \mathcal{L}(H_1,Y)$ is an admissible output operator for $\mathbb T$, there exists an extension given by:

\begin{equation}
\label{eq:extensionC}
    C_\Lambda z:= \lim_{\lambda\rightarrow + \infty} C\lambda (\lambda \mathrm I_H-A)^{-1}z,
\end{equation}
with $\lambda$ a real number and $z\in D(C_\Lambda)$, where the space $D(C_\Lambda)$ is defined as follows
\begin{equation}
    D(C_\lambda):=\lbrace z\in H\mid \text{the limit \eqref{eq:extensionC} exists}\rbrace.
\end{equation}
Note that it is not the only extension of $C$, but it turns out that this one is sufficient for our purpose. For other extensions, we refer the interested reader to \cite{weiss1994transfer}. Let us now give a definition of what is a well-posed LTI system (see e.g., \cite{tucsnak2014well} for a general introduction of such a notion):

\begin{definition}[Well-posed LTI systems]
\label{def:well-posed}
A system $\Sigma$ defined as in \eqref{eq:ISS} is called well-posed if for some $t>0$ there exists a positive constant $M_t$ such that, for almost every $t$:
$$
\Vert z(t)\Vert_H + \Vert y\Vert_{L^2([0,t];Y)} \leq M_t (\Vert z_0\Vert_H + \Vert u\Vert_{L^2([0,t];U)}).
$$
\end{definition}

It is well known that a system without output (or with a bounded output operator) is well posed as soon as $B$ is an admissible operator for the $\mathbb T$ in the sense given in \cite{tucsnak2009observation}. It is also the case for system without input and just an output: the system is well-posed if the output operator is admissible for $\mathbb T$. However, when output and input are involved, another property is needed to ensure the well-posedness of \eqref{eq:ISS}. Among the sufficient conditions for such a well-posedness, we may mention that a triple $(A,B,C)$ is regular, for which we give now a definition. 
\begin{definition}
\label{ass:structure}
The triple $(A,B,C)$ is called regular if the following properties are satisfied:
\begin{itemize}
    \item[1.] The operator $B$ is admissible for the semigroup $(\mathbb T_t)_{t\geq 0}$. 
    \item[2.] The operator $C$ is admissible for the semigroup $(\mathbb T_t)_{t\geq 0}$.
    \item[3.]  The function $\mathbf{G}(s) = C_\Lambda(s\mathrm I_H-A)^{-1} B$ makes sense for some (hence for every) $s\in \rho(A)$ (where $\rho(A)$ denotes the resolvent set of $A$).
    \item[4.] The function $s\mapsto \Vert \mathbf G(s)\Vert_{\mathcal{L}(U,Y)}$ is bounded on some right half plane. 
\end{itemize}
\end{definition}
An equivalent characterization (in particular, Item 4.) of a well-posed system that is regular is that the following limit $\lim_{\lambda \rightarrow + \infty} \mathbf{G}(\lambda)u$ exists for $\lambda$ a real number and for every $u\in U$. This limit describes the feedthrough operator $D\in \mathcal L(U,Y)$, i.e.:
\begin{equation}
Du:=\lim_{\lambda\rightarrow + \infty} G(\lambda) u,\: \forall u\in U,
\end{equation}
therefore the transfer function of a regular system is given by:
\begin{equation}
\mathbf H(s) = C_\Lambda(s\mathrm I_H - A)^{-1} + D.
\end{equation}
This means that from the triple $(A,B,C)$, one can deduce the operator $D$. Therefore, the characterization of a system (that is infinite-dimensional or not) relies on the triple $(A,B,C)$. Furthermore, it implies that, even if \eqref{eq:ISS} does not admit any feedthrough operator $D$, extending $C$ may make appear such an operator.

We now give some notation. We call $\mathbf P_\tau$ the truncation operator applied to any signals $u\in L^2([0,\infty])$. It is defined as follows:
\begin{equation*}
\mathbf P_\tau u= \left\{
\begin{aligned}
&u(t),\: &\forall &t\in [0,\tau]\\
&0, &\forall &t >\tau
\end{aligned}
\right.
\end{equation*}
We will denote the following spaces as follows $\mathscr{U}:=L^2_{\mathrm{loc}}([0,+\infty);U)$ and $\mathscr{Y}:=L^2_{\mathrm loc}([0,+\infty);U)$ and consider that $u\in \mathscr{U}$ and $y\in \mathscr{Y}$. It follows that $$\mathbf P_\tau u\in L^2([0,\infty);U),$$ since it vanishes for all $t>\tau$. Assuming that $(A,B,C)$ is regular, we now define a system such as the one in \eqref{eq:ISS} as follows:

\begin{equation}
\begin{bmatrix}
z(\tau) \\ \mathbf P_\tau y
\end{bmatrix}= \begin{bmatrix}
\mathbb T_\tau & \mathbb \phi_\tau\\
\psi_\tau & \mathbb F_\tau
\end{bmatrix}\begin{bmatrix}
z_0\\
\mathbf P_\tau u
\end{bmatrix},
\end{equation}
where $(\mathbb T_t)_{t\geq 0}\in \mathcal L(H)$ is the semigroup generated by $A$, $(\phi_t)_{t\geq 0}\in \mathcal L(\mathscr U,H)$ is given by
\begin{equation}
\label{eq:phi}
    \phi_\tau \mathbf P_\tau u=\int_0^\tau \mathbb T_{t-s} B u(s) \diff s,
\end{equation}
the operator $(\psi_t)_{t\geq 0}\in\mathcal L(H,\mathscr Y)$ is given by
\begin{equation}
\label{eq:psi}
    \psi_\tau z_0 := \int_0^\tau C_\Lambda \mathbb T_s z_0 \diff s,
\end{equation}
and the operator $(\mathbb F_t)_{t\geq 0}\in\mathcal{L}(\mathscr U,\mathscr Y)$ is defined as follows:
\begin{equation}
\label{eq:F}
\mathbb F_\tau \mathbf P_\tau u = C_\Lambda \int_0^\tau \mathbb T_{\tau-s} B \mathbf P_\tau u(s) \diff s + D\mathbf P_\tau u.
\end{equation}
We will call further $\mathbb F_\infty$ the operator mapping $L^2_{\mathrm loc} ([0,\infty);U)$ to $L^2_{\mathrm loc} ([0,\infty);Y)$ and given by:

\begin{equation}
\label{eq:Finf}
\mathbb F_\infty u =  C_\Lambda \int_0^\tau \mathbb T_{\tau-s} B u(s) \diff s + D u,
\end{equation}
which is valid for almost every $t\geq 0$ and for every $u\in L^2_{\mathrm loc}([0,\infty);U)$.

The boundedness property of these operators is deduced from the assumption that the triple $(A,B,C)$ is well-posed as in Definition \ref{ass:structure}. The fact that the operators $\mathbb \psi$ and $\mathbb F$ can be explicitly related to the operators $(A,B,C)$ is an implication of Item 3. of Definition \ref{ass:structure} (namely, the fact that the triple $(A,B,C)$ is regular). 

Now, we need some concepts of solutions, which are not easy in general to consider when adding inputs and outputs. Let us define, therefore, the space $D(S)$, which is densely defined in the space $H\times U$, and which is equipped with the following graph norm norm:
$$
\Vert (z,u)\Vert^2_{D(S)}:= \Vert z\Vert_H^2 + \Vert u\Vert_U^2 + \Vert Az+Bu\Vert_H^2,
$$
where the related operator $S$ can be defined as:
$$
S\begin{bmatrix}
    z\\ u
\end{bmatrix}:=\begin{bmatrix}
    A & B\\
    C & 0
\end{bmatrix}\begin{bmatrix}
    z\\ u
\end{bmatrix}.
$$
Therefore, the system \eqref{eq:ISS} can be written as follows:
\begin{equation}
    \begin{bmatrix}
        \frac{\diff}{\diff t} z(t)\\
        y(t)
    \end{bmatrix} = S \begin{bmatrix}
        z(t)\\
        u(t)
    \end{bmatrix}.
\end{equation}
\begin{definition}
\label{def:solutions}
The triple $(z,u,y)$ is called a strong solution to \eqref{eq:ISS} if:
\begin{itemize}
    \item[1.] $z\in C^1([0,\infty);H)$.
    \item[2.] $u\in C([0,\infty);H)$, $y\in C([0,\infty);Y)$.
    \item[3.] $\begin{bmatrix}
        z(t) \\ u(t)
    \end{bmatrix}\in D(S)$ for all $t\geq 0$. 
\end{itemize}
The triple $(z,u,y)$ is called a generalized solution to \eqref{eq:ISS} if:
\begin{itemize}
    \item[1.] $z\in C([0,\infty);H)$.
    \item[2.]  $u\in L^2_{\mathrm loc}([0,\infty);U)$, $y\in L^2_{\mathrm loc}([0,\infty);Y)$.
    \item[3.] There exists a sequence $(z_k,u_k,y_k)$ of classical solutions to \eqref{eq:ISS} such that $z_k\rightarrow z$ in $C([0,\infty);H)$, $u_k\rightarrow u$ in $L^2_{\mathrm loc}([0,\infty);U)$ and $y_k\rightarrow y$ in $L^2_{\mathrm loc}([0,\infty);Y)$.  
\end{itemize}
\end{definition}

The existence of strong solution can be stated as soon as $\begin{bmatrix}z_0\\ u(0)\end{bmatrix}\in D(S)$ as stated in the following result, taken from \cite[Proposition 4.3.]{tucsnak2014well}:
\begin{proposition}
If $u\in C^2([0,\infty);U)$ and $\begin{bmatrix}z_0\\ u(0)\end{bmatrix}\in D(S)$, then the system \eqref{eq:ISS} admits a unique strong solution satisfying $z(0) = z_0$. 
\end{proposition}

We have also existence of generalized solutions in the case of regular systems (as stated in \cite[Theorem 5.6]{tucsnak2014well}):
\begin{proposition}
For any $u\in L^2_{loc}([0,\infty);U)$ and any $z_0\in H$, there exists a unique generalized solution to \eqref{eq:ISS} $\begin{bmatrix}
    z(t)\\
    y(t)
\end{bmatrix}\in H\times L^2_{\mathrm loc}([0,\infty);Y)$ for almost every $t\geq 0$. This generalized solution is expressed through the semigroup $\mathbb T$ generated by $A$ and the operators given in \eqref{eq:phi}, \eqref{eq:psi} and \eqref{eq:F}. 
\end{proposition}

 For further details about well-posed infinite-dimensional systems, we refer the interested reader to \cite{tucsnak2014well,staffans2002passive,weiss1994transfer}. 

 \begin{example}[Korteweg-de Vries equation]
 \label{example:kdv}
Consider the following Korteweg-de Vries equation:
\begin{equation}
    \left\{
\begin{aligned}
& z_t + z_x + z_{xxx} = 0,\: &\text{ on } &[0,\infty)\times [0,L],\\
& z(t,0) =z(t,L)= 0, \: &\text{ on } &[0,\infty),\\
& z_x(t,L) = u(t), &\text{ on } & [0,\infty),\\
& y(t) = z_x(t,0), &\text{ on } & [0,\infty),\\
& z(0,x) = z_0(x), &\text{ on } & [0,L]. 
\end{aligned}
    \right.
\end{equation}
The state space $H$ is given by $H:=L^2(0,L)$, the control space and the output space are $U=Y=\mathbb R$, and the operator $A$ is defined as:
$$
Az = - z^\prime - z^{\prime\prime\prime},
$$
with the domain:
$$
D(A):=\lbrace z\in H^3(0,L)\mid z(0) = z(L)= z^\prime(L) = 0\rbrace, 
$$
which can be proved to be the generator of a strongly continuous semigroup of contractions (see e.g., \cite{cerpa2013control,rosier1997kdv,nguyen2024rapid}). The operator $C$ is defined as 
$$
C z = \delta^\prime_{x=0} z=z^\prime(0),
$$
where $\delta^\prime_{x=0}$ is the evaluation operator of the derivative of $z$ at $x=0$. The control operator $B$ can be defined as the Delta operator given by $\langle Bu,w\rangle_{D(A^*)',D(A^*)}=u \cdot w^\prime(L)$, where $A^*$ denotes the adjoint operator of $A$ ad $D(A^*)$ its domain, that are defined as follows:
$$
A^* w=w^\prime + w^{\prime\prime\prime},\: D(A^*):=\lbrace w\in H^3(0,L)\mid w(0)=w(L)=w^\prime(0)=0\rbrace.
$$
Note that $D(A^*)'$, the dual of $D(A^*)$, can be identified by $H_{-1}$. 

As soon as $L\notin \mathcal N$, where $\mathcal N$ is given by:
\begin{equation}
    \mathcal N:=\left\{2\pi\sqrt{\tfrac{k^2+kl+l^2}{3}}:k,l\in\mathbb N\right\},
\end{equation}
we know that the output $y$ never vanishes (see e.g., \cite{rosier1997kdv}). Moreover, the time derivative of the energy $E(z):=\frac{1}{2} \Vert z\Vert_{L^2(0,L)}$ satisfies
$$
\frac{\diff}{\diff t} E(z) = |u(t)|^2 - |y(t)|^2,
$$
which proves that the system is scattering passive. Therefore, it is well-posed (see e.g., \cite{tucsnak2014well,staffans2005well,singh2024local} for an introduction on this topic). Indeed, taking the integral of the latter equation, one obtains that, for all $t\geq 0$:
\begin{equation}
\label{eq:kdv-scattering}
\Vert z(t)\Vert_{L^2(0,L)}^2 - \Vert z_0\Vert_{L^2(0,L)}^2 = \int_0^t |u(s)|^2\diff s - \int_0^t |y(s)|^2\diff s.
\end{equation}
Therefore, the inequality given in Definition \ref{def:well-posed} is satisfied, meaning in particular that $B$ and $C$ are admissible for the semigroup generated by $A$.

To compute the transfer function of the KdV equation, let us take a look at the following boundary value problem:
\begin{equation}
\label{eq:Laplace-KdV}
\left\{
\begin{aligned}
    & s \hat z(s,x) + \hat z_x(s,x) + \hat z_{xxx}(s,x) = 0,\\
    & \hat z(s,0) = \hat z(s,L) = 0,\\
    & \hat z_x(s,L) = \hat u(s),\\
    & \hat z(0,x) = 0,
\end{aligned}
\right.
\end{equation}
where $\hat z$ denotes the Laplace transform of $z$ with respect to $t$ and $\hat u$ denotes the Laplace transform of $u$ and where $s \in \mathbb C$. Suppose that $s\in \rho(A)$ with $\rho(A)$ the resolvent set of $A$, so that one can invert the operator $s\mathrm I_H - A$ (which corresponds to the resolvent operator).

The characteristic equation
$$
s + r^{3} + r = 0,
$$
admits three roots denoted by $r_1,r_2,r_3\in \mathbb C$. We will see that the exact formula of $r_1,r_2,r_3$ is not needed for our purpose. The solution to \eqref{eq:Laplace-KdV} can be written as:
$$
z(x) = c_1 e^{r_1 x} + c_2 e^{r_2 x} + c_3 e^{r_3 x},
$$
where $c_1,c_2,c_3\in \mathbb C$ are constants to be determined with the boundary conditions given in \eqref{eq:Laplace-KdV}. The coefficients satisfy the following algebraic equation:
\begin{equation}
\begin{bmatrix}
    1 & 1 & 1\\
    e^{r_1 L} & e^{r_2 L} & e^{r_3L}\\
    r_1 e^{r_1L} & r_2 e^{r_2 L} & r_3 e^{r_3L}
\end{bmatrix}\begin{bmatrix}
    c_1 \\ c_2 \\ c_3
\end{bmatrix} = \begin{bmatrix}
0 \\ 0 \\ \hat u(s)
\end{bmatrix}. 
\end{equation}
Using the Cramer's rule, one obtains the following result:
\begin{equation}
\begin{aligned}
&c_1 := \hat u(s) \frac{e^{r_3 L}-e^{r_2 L}}{(e^{r_3L}-e^{r_2L})r_1 + (e^{r_1L} - e^{r_3L})r_2 + (e^{r_2L}-e^{r_1L})r_3},\\
&c_2 := \hat u(s) \frac{e^{r_1L}-e^{r_3L}}{(e^{r_3L}-e^{r_2L})r_1 + (e^{r_1L} - e^{r_3L})r_2 + (e^{r_2L}-e^{r_1L})r_3},\\
& c_3 = \hat u(s) \frac{e^{r_2 L}-e^{r_1 L}}{(e^{r_3L}-e^{r_2L})r_1 + (e^{r_1L} - e^{r_3L})r_2 + (e^{r_2L}-e^{r_1L})r_3}.
\end{aligned}
\end{equation}
Moreover, the Laplace transform of the output is given by
$$
\hat y (s) = z^\prime(s,0) = c_1 r_1 + c_2 r_2 + c_3 r_3.
$$
Surprisingly, it turns out that the transfer function, given by $\mathbf H(s) = \frac{\hat y(s)}{\hat u(s)}=1$ (it is therefore, a constant transfer function), which shows that, for all $s$ real, there exists a strong limit for $\mathbf H(s)$ as $s\rightarrow + \infty$. This strong limit is $1$ and is equal to $D$.  
 \end{example}

\subsection{Building a coercive Lyapunov functional}

Before considering the Lyapunov functional, let us give some notation. As soon as we consider the time derivative of the Lyapunov functionals (for example, $V:H\rightarrow \mathbb R_+$), we will consider its Dini derivative, given by:
$$
\frac{\diff}{\diff t}V(z):=\lim_{t \rightarrow 0^+} \frac{1}{t} V(z(t))-V(z_0),
$$
where $z(t)$ corresponds, for every $t\geq 0$, to a trajectory of a well-posed system.

As we are going to build a coercive ISS Lyapunov functional, we focus now on some details about Lyapunov functionals for infinite-dimensional systems. It is known since decades (see e.g., \cite{datko1970extending}) that as soon as $(\mathbb T)_{t\geq 0}$ is an exponentially stable semigroup, then there exists a self-adjoint operator $P_0\in\mathcal L(H)$ and a positive constant $\alpha$ such that, for all $z\in D(A)$
\begin{equation}
\label{eq:lyap}
0<\langle P_0 z,z\rangle_H\leq \Vert P_0\Vert_{\mathcal L(H)} \Vert z\Vert_H^2,\quad 2 \langle P_0z,Az\rangle_H \leq -\alpha \Vert z\Vert^2_H.
\end{equation}
Now, looking at the following system:
\begin{equation}
\label{eq:dyn-input}
\left\{
\begin{aligned}
& \frac{\diff}{\diff t} z(t) = Az(t) + Bu(t),\: t\in [0,+\infty)\\
& z(0)=z_0,
\end{aligned}
\right.
\end{equation}
then it turns out that if $PA$ admits an extension called to a bounded operator in $H$ (called $PA$ again), then the Lyapunov functional defined in \eqref{eq:lyap} is ISS for \eqref{eq:dyn-input}, i.e. there exists $\alpha_1,\alpha_2>0$ such that, for all $z\in D(A)$
\begin{equation}
\label{iss:jacob}
\frac{\diff}{\diff t} W(z) \leq - \alpha_1 \Vert z\Vert_H^2 + \alpha_2 \Vert u\Vert^2_{U}, 
\end{equation}
with $W(z) =:\langle Pz,z\rangle_H$. This result is given in \cite[Theorem 5.3.]{jacob2020noncoercive}. We now provide a result allowing to construct a coercive ISS Lyapunov functional for \eqref{eq:ISS} by assuming that the pair $(A,C)$ is exactly observable.
\begin{theorem}
\label{thm:ISS-coercive}
Suppose that the triple $(A,B,C)$ is regular and, moreover, suppose that the pair $(A,C)$ is exactly observable. Then, there exist $\alpha_1,\alpha_2,\alpha_3$ and a coercive operator $P\in \mathcal L(H)$\footnote{Note that the coercivity property means that there exists $\beta>0$ such that $\langle Pz,z\rangle\geq \beta \Vert z\Vert_H^2$ for any $z\in H$.} such that for any strong solutions to \eqref{eq:ISS}:
\begin{equation}
\label{eq:coercive-ISS}
\frac{\diff }{\diff t} V(z)\leq - \alpha_1 \Vert z\Vert^2_H + \alpha_2 \Vert u\Vert^2_U - \alpha_3 \Vert Cz\Vert_Y^2. 
\end{equation}
Any generalized solutions to \eqref{eq:ISS} satisfy, for every $z_0\in H$, every $u\in L^2_{\mathrm loc}([0,\infty);U)$:
\begin{equation}
\label{eq:lyap-generalized}
    V(z(t)) \leq e^{- \alpha_1t} V(z_0) + \alpha_2 \int_0^t e^{-\alpha(t-s)} \Vert u(s)\Vert_U^2 \diff s - \alpha_3 \int_0^t \Vert Cz(s)\Vert_Y \diff s,
\end{equation}
where $\alpha_1$, $\alpha_2$, $\alpha_3$ are positive constants. 
\end{theorem}

\begin{remark}
It is important to note that one only needs the existence of an operator $C$ such that the pair $(A,C)$ is exactly observable. In other words, bounded operators $C$ are also good candidates, which implies that the assumption is not that strong. However, in some cases, such operators do not exist, such as the parabolic case. However, this specific case has been tackled in \cite{mironchenko2023coercive}, and it turns out that the usual norm for analytic systems is a coercive ISS Lyapunov functional if the associated semigroup generated by $A$ is supposed to be contractive.
\end{remark}

\begin{proof}
All along this proof, we will consider strong solutions, as the ones introduced in Definition \ref{def:solutions}, i.e., we have $(z_0,u(0))\in D(S)$ and $u\in C^2([0,\infty);U)$. The result about generalized solution follows from a standard density argument. We will first follow the idea given in \cite{lorenzetti2023saturating}, upon which is based the notion of observability Gramian \cite[p. 140]{tucsnak2009observation}. To do so, consider the operator $P_1\in\mathcal L(H)$ defined as:
\begin{equation}
\label{def:P1}
P_1z:=\lim_{\tau \rightarrow + \infty} \int_0^\tau \mathbb T^*_t C^*C\mathbb T_t z \diff t.
\end{equation}
We also define $P_1^\tau:=\int_0^\tau \mathbb T^*_t C^*C\mathbb T_t z \diff t.$ Using the admissibility of the operator $C$, one obtains that, for every $\tau >0$ and every $z\in H$:
$$
\langle P^\tau_1 z,z\rangle_H =  \int_0^\tau \Vert C \mathbb T_t z\Vert^2_U \diff t \leq K_C \Vert z\Vert_H^2, 
$$
where $K_C$ is independent on $\tau$ since $A$ is exponentially stable (see e.g., \cite[Proposition 4.4.5]{tucsnak2009observation}). This means that, taking the limit when $\tau$ goes to infinity, one obtains the following upper-bound:
$$
\langle P_1 z,z\rangle_H \leq K_C \Vert z\Vert^2_H.
$$
Using now the exact observability of the operator with respect to the semigroup $\mathbb T$ (\cite[Definition 6.5.1]{tucsnak2009observation}), we know that there exists $\tilde K_C>0$ such that
$$
\langle P_1 z,z\rangle_H =  \lim_{\tau \rightarrow + \infty}\int_0^\tau \Vert C \mathbb T_t z\Vert^2_Y \diff t \geq \tilde K_C \Vert z\Vert^2_H,
$$
implying therefore that the functional $\langle P_1z,z\rangle_H$ is coercive.
According to \cite[Theorem 5.1.1.]{tucsnak2009observation}, since $C$ is infinite-time admissible for the semigroup $\mathbb T$, the operator $P_1$ satisfies the following inequality
\begin{equation}
\label{eq:lyapgram}
2\langle P_1z,Az\rangle\leq - \Vert Cz\Vert^2_Y,\: \forall z\in D(A).
\end{equation}
It means therefore that the Lyapunov functional:
$$
V(z):=\langle Pz,z\rangle,
$$
with $P=P_0+P_1$ is a coercive Lyapunov functional for \eqref{eq:dyn-input}, where $P_0\in \mathcal L(H)$ is given in \eqref{eq:lyap}. From \cite[Theorem 5.3.]{jacob2020noncoercive}, we already know that the inequality \eqref{iss:jacob} holds. Let us now check whether such an inequality holds also for the operator $P_1$. We first note that the solution to \eqref{eq:ISS} is expressed as:
$$
z(t) = \mathbb T_t z_0 + \int_0^t \mathbb T_{t-s} B \mathbf P_t u(s) \diff s = \mathbb T_t z_0 + \phi_t u,
$$
where $\phi_t$ is defined in \eqref{eq:phi}. 
Pick any $u\in \mathscr U$ and $T>0$ arbitrary. Since we consider initial condition $z_0\in D(A)$, then $C_\Lambda = C$, meaning that $y(t)= C z(t) + Du(t)$. Then, one has, for all $t\geq 0$ and for all $(z_0,u(0))\in D(S)$
\begin{align*}
\langle P_1^T z(t),z(t)\rangle_H - \langle P_1^T z_0,z_0\rangle_H = &\int_0^T \left\Vert C \mathbb T_\tau z(t) \right\Vert^2_{Y}\diff \tau - \int_0^T \Vert C \mathbb T_\tau z_0 \Vert^2_{Y} \diff \tau\\
= & \int_0^T \left \Vert C \mathbb T_{\tau + t} z_0 + C \mathbb T_\tau \int_0^t \mathbb T_{(t-s)} B\mathbf P_t u(s) \diff s \right\Vert_{Y}^2 \diff \tau \\
&- \int_0^T \Vert C \mathbb T_\tau z_0 \Vert^2_{Y} \diff \tau\\
= & \int_0^T \Vert C \mathbb T_{\tau + t} z_0\Vert^2_{ Y} \diff \tau \\
&+ \int_0^T \left \Vert C \mathbb T_\tau \int_0^t \mathbb T_{(t-s)} B\mathbf P_t u(s) \diff s \right\Vert_{ Y}^2 \diff \tau \\
&  + 2 \int_0^T \left\langle  C \mathbb T_{\tau + t} z_0, C \mathbb T_\tau \int_0^t \mathbb T_{(t-s)} B\mathbf P_t u(s) \diff s\right \rangle_{Y}\\
&- \int_0^T \Vert C \mathbb T_\tau z_0 \Vert^2_{\mathscr Y} \diff \tau\\
\leq & \int_0^T \Vert C \mathbb T_{\tau + t} z_0\Vert^2_{Y} \diff \tau \\
&+ \int_0^T \left \Vert C \mathbb T_\tau \int_0^t \mathbb T_{(t-s)} B\mathbf P_t u(s) \diff s \right\Vert_{Y}^2 \diff \tau \\
& + 2 \int_0^T \Vert C \mathbb T_{\tau + t} z_0 \Vert_{Y} \left\Vert  C \mathbb T_\tau \int_0^t \mathbb T_{(t-s)} B\mathbf P_t u(s) \diff s \right \Vert_{Y} \diff \tau \\
&- \int_0^T \Vert C \mathbb T_\tau z_0 \Vert^2_{Y} \diff \tau,
\end{align*}
where we have used the Cauchy-Schwarz inequality to obtain the last inequality.

We know that the operator $\int_0^T \Vert C \mathbb T_\tau z\Vert_H^2 \diff \tau$ is uniformly bounded (i.e., this bound does not depend on $T$). The terms involving the input $u$ are also uniformly bounded by $T$ by taking $z = \int_0^t \mathbb T(t-s) u(s) \diff s$, which belongs to $H$ (since $B$ is an admissible operator for $\mathbb T$). Therefore, one can take $T$ going to infinity, and then:

\begin{align*}
\langle P_1 z(t),z(t)\rangle_H - \langle P_1 z_0,z_0\rangle_H \leq & \int_0^\infty \Vert C \mathbb T_{\tau + t} z_0\Vert^2_{ Y} \diff \tau \\
&+ \int_0^\infty \left \Vert C \mathbb T_\tau \int_0^t \mathbb T_{(t-s)} B\mathbf P_t u(s) \diff s \right\Vert_{ Y}^2 \diff \tau \\
& + 2 \int_0^\infty \Vert C \mathbb T_{\tau + t} z_0 \Vert_{Y} \left\Vert  C \mathbb T_\tau \int_0^t \mathbb T_{(t-s)} B\mathbf P_t u(s) \diff s \right \Vert_{Y} \diff \tau \\
&- \int_0^\infty \Vert C \mathbb T_\tau z_0 \Vert^2_{Y} \diff \tau,\\
 \leq &\int_0^\infty \Vert C \mathbb T_{\tau + t} z_0\Vert^2_{ Y} \diff \tau \\
 &+ \left(1 + \frac{1}{\delta}\right) \int_0^\infty \left \Vert C \mathbb T_\tau \int_0^t \mathbb T_{(t-s)} B\mathbf P_t u(s) \diff s \right\Vert_{Y}^2 \diff \tau\\
& + \delta\int_0^\infty \Vert C T_{\tau + t} z_0\Vert_Y^2 \diff \tau - \int_0^\infty \Vert C \mathbb T_\tau z_0 \Vert^2_{Y} \diff \tau,
\end{align*}
where we have used the Young's inequality in the last inequality. The constant $\delta$ will be chosen later on. 

Note that the operator $C \mathbb T_\tau \int_0^t \mathbb T_{(t-s)} Bu(s) \diff s$ is related with $\mathbb F$, because
$$
 \int_0^\infty C \mathbb T_\tau \int_0^t \mathbb T_{(t-s)} Bu(s) \diff s \diff \tau =\mathbb F_\infty u - Du,
$$
where $\mathbb F_\infty$ is defined in \eqref{eq:F} and \eqref{eq:Finf}. Since the triple $(A,B,C)$ is regular, and because $\mathbb T$ is an exponentially stable semigroup, then according to \cite[Proposition 4.1.]{weiss1994transfer}, there exist $M_\gamma\geq 1$ such that
$$
\Vert \mathbb F_\infty u\Vert_\mathscr Y \leq M_\gamma \Vert \phi_t u\Vert_H, 
$$
where $\phi$ is defined in \eqref{eq:phi}. Since $B$ is admissible for $\mathbb T$, and because $\mathbb T$ is an exponentially stable semigroup, there exists $K_B>0$ independent on $t$ such that
\begin{equation}
\Vert \mathbb F_\infty u\Vert_\mathscr Y \leq K_B M_\gamma   \int_0^\infty  \Vert \mathbf P_t u\Vert_U \diff s
\end{equation}
Since $D$ is a bounded operator $D\in \mathcal L(U,Y)$, one can easily prove that:
$$
\Vert Du\Vert_{Y} \leq \Vert D\Vert_{\mathcal L(U,Y)}  \int_0^t  \Vert \mathbf P_t u\Vert_U^2 \diff t.
$$
Moreover, $\mathbb T$ is an exponential semigroup, then there exists $M$ such that
$$
\Vert C \mathbb T_{t} z_0 \Vert_{\mathscr Y} \leq M \Vert C z_0\Vert_\mathscr Y = M\int_0^t \mathbf P_t C \mathbb T_\tau z_0 \diff \tau 
$$

Finally, one can prove that, setting $K:= \left(1+\frac{1}{\delta}\right)(\Vert D\Vert_{\mathcal L(U,Y)} + K_B M_\gamma)$, for all $z_0\in D(A^2)$:
\begin{align*}
\langle P_1 z(t),z(t)\rangle_H - \langle P_1 z_0,z_0\rangle_H \leq & \Vert C \mathbb T_{t} z_0\Vert^2_{\mathscr Y} - \Vert C z_0\Vert^2_{\mathscr Y} + \delta M \int_0^t \Vert \mathbf P_t \psi_t z_0\Vert^2_Y \diff t\\
& + K \int_0^t \Vert \mathbf P_t u\Vert_U \diff t
\end{align*}

Now we are going to divide the latter inequality by $t$ and take the limit. First, note that:
$$
\lim_{t\rightarrow 0} \frac{\langle P_1 z(t),z(t)\rangle_H - \langle P_1 z_0,z_0\rangle_H}{t} = \frac{\diff}{\diff t} U(z_0),
$$
where $U(z):=\langle P_1 z,z\rangle_H$ (defined for any $z\in H$) and where the differential is understood as a Dini derivative. Second, by definition of the generator of a semigroup, one obtains the following limit:
$$
\lim_{t\rightarrow 0}\frac{\Vert C \mathbb T_{t} z_0\Vert^2_{\mathscr Y} - \Vert C z_0\Vert^2_{\mathscr Y}}{t} = 2\langle P_1 A z_0,z_0\rangle_H \leq  - \Vert C z_0\Vert^2_H,
$$
where the last equality comes from \eqref{eq:lyapgram}. Moreover, since $u\in \mathscr U$, the input $u$ admits a Lebesgue point at $0$ (as stated in \cite[Lemma 5.4.10]{staffans2005well}), meaning that, for almost every $t\geq 0$:
$$
 K \lim_{t\rightarrow 0}\int_0^t  \Vert \mathbf P_t u\Vert_U^2 \diff s = K\Vert u(t)\Vert_U.
$$
Moreover, since the triple $(A,B,C)$ is regular, then the following limit exists:
$$
\delta M \lim_{t\rightarrow 0} \frac{1}{t} \int_0^t \Vert \mathbf P_t \psi_t z_0\Vert_Y \diff t = \delta M\Vert C z_0\Vert^2_Y
$$
Then, gathering all these limits, one obtains the inequality, satisfied for every $z_0$ and $u(0)$ such that $Az_0 + Bu(0)\in H$, every $u\in L^2_{\mathrm loc}([0,\infty[;U)$ and for almost every $t\geq 0$
\begin{equation}
    \frac{\diff }{\diff t} U(z_0) \leq - (1-\delta M) \Vert Cz_0\Vert^2_Y + K \Vert u(t)\Vert^2_U,
\end{equation}
To conclude, it suffices to take $\delta<\frac{1}{M}$. Therefore, the Lyapunov functional $V(z) = U(z) + W(z)$, where $W$ is given in \eqref{eq:lyap}, satisfies the inequality \eqref{eq:coercive-ISS} for all $z_0\in D(S)$ with $u\in C^2([0,\infty;U)$. By a density argument, one can show that the inequality \eqref{eq:lyap-generalized} for every generalized solution (i.e., for all $z_0\in H$ and $u\in \mathscr U$). 
\end{proof}

As a corollary (with a direct proof following the arguments below and the ones provided in \cite{jacob2020noncoercive}), one can prove the following:

\begin{corollary}
\label{corollary-strict}
Consider that $u=0$ in \eqref{eq:ISS}. Suppose that there exist two functionals $$V_1(z):=\langle P_1z,z\rangle_H,V_0(z):=\langle P_0 z,z\rangle_H,$$ with $P_i\in\mathcal L(H)$ that are two self-adjoint operators with $i\in \lbrace 0,1\rbrace$, mapping $H$ to $\mathbb R_+$, and which satisfy, along the trajectories of \eqref{eq:ISS} (with $u=0$):
$$
\frac{\diff}{\diff t} V_1(z) \leq - \alpha_2 \Vert Cz \Vert^2_Y,\: \frac{\diff}{ \diff t} V_2(z) \leq - \alpha_2 \Vert z\Vert^2_H,
$$
where $\alpha_1,\alpha_2$ are positive constants. We consider that both operators $P_1,P_2$ admit extensions for $P_1A,P_2A$ (still called $P_1A$, $P_2A$). 

Now, consider that $u\neq 0$, with $u\in \mathscr U$. Then, there exists $\tilde \alpha_1,\tilde \alpha_2$ positive constants such that the Lyapunov function $V:=V_1+V_2$ satisfies along the trajectories of \eqref{eq:ISS} with $u\in\mathscr U$:
$$
\frac{\diff}{\diff t} V(z) \leq - \tilde \alpha_1 \Vert z\Vert_H^2 - \tilde\alpha_2 \Vert Cz\Vert^2_H + \alpha_3 \Vert u(t)\Vert^2_U,
$$
for almost every $t\geq 0$. If the pair $(A,C)$ is exactly observable, the Lyapunov functional $V$ is coercive. 
\end{corollary}

\begin{remark}
Note that in some cases one only has the Lyapunov functional given by $V_1$ in the statement of the Corollary \ref{corollary-strict}. It is for instance the case of scattering passive system, such as the KdV equation given in \eqref{example:kdv}. Such storage functions are called "weak Lyapunov functional" in the literature, and there exist strictification techniques in order to make these weak Lyapunov functionals strict, i.e. a Lyapunov functional satisfying the statement of Theorem \ref{thm:ISS-coercive}. One of them is written in \cite{praly2019observers}, and this technique, usually used in the finite dimensional context, has been extended to the case of the KdV equation in \cite{balogoun2021iss}. This strictification technique is based on an observer, implying therefore that an observability property has to be satisfied for the system for which we aim at building a strict Lyapunov functional.
\end{remark}

\section{Applications}

\label{sec:applications}

\subsection{Singular perturbation}

We are going to study the following abstract system:
\begin{equation}
\label{eq:spm}
    \left\{
\begin{aligned}
     &\varepsilon \frac{\diff}{\diff t}z = A_1 z + B_1C_2 w,\\
     & \frac{\diff}{\diff t} w = A_2 w + B_2 C_1 z,\\
     & y = \begin{bmatrix}
         C_1 z & C_2 w
     \end{bmatrix},\\
     & z(0) = z_0,\: w(0) = w_0.
\end{aligned}
    \right.
\end{equation}
where $A_1:\: D(A_1)\subset H\rightarrow H$, $B_1\in \mathcal{L}(\mathbb R^p,H_{-1})$, $C_1\in \mathcal{L}(D(A),\mathbb R^m)$, $A_2\in \mathbb R^{n\times n}$, $B_2\in \mathbb R^{m\times n}$, $C_2\in \mathbb R^{n\times p}$ and $\varepsilon>0$. Therefore, it is a singularly perturbed system composed by an infinite-dimensional system (described with the state variable $z$) and a finite-dimensional system (described with the state variable $w$). Rewriting \eqref{eq:spm} as follows:

\begin{equation}
\left\{
\begin{aligned}
& \frac{\diff}{\diff t} \begin{bmatrix} z \\ w\end{bmatrix} = \begin{bmatrix}
\frac{1}{\varepsilon} A_1 & \frac{1}{\varepsilon} B_1 C_2\\
B_2C_1 & A_2
\end{bmatrix}\begin{bmatrix}
    z \\ w
\end{bmatrix},\\
& y = \begin{bmatrix}
         C_1 z & C_2 w
     \end{bmatrix},
\end{aligned}
\right.
\end{equation}

We suppose that this system is well-posed and regular (conditions for such a well-posedness are given, for instance, in \cite{weiss2009well,weiss1994regular,weiss1997dynamic}). In particular, we state the following assumption:
\begin{assumption}
The following properties are satisfied:
\begin{itemize}
    \item[1.] $A_1$ generates a strongly continuous semigroup $(\mathbb T_t)_{t\geq 0}$.
    \item[2.] The operator $B_1$ is infinite-time admissible for the semigroup $\mathbb T$.
    \item[3.] The operator $C_1$ is admissible for the semigroup $\mathbb T$.
    \item[4.] The triple $(A_1,B_1,C_1)$ is regular.
    \item[5.] Consider the transfer function:
    $\mathbf H_1(s):= C_{\Lambda,1}(s\mathrm I_{\mathbb R^{n\times n}}-A_1)^{-1}B_1 + D_1$, with $$D_1:=\lim_{\lambda \rightarrow + \infty} C_{\Lambda,1}(\lambda\mathrm I_{\mathbb R^{n\times n}}-A_1)^{-1}B_1,$$ where $C_{\Lambda,1}$ is the extension of $C_1$ and $\lambda$ a real number. Consider also $\mathbf H_2(s):= C_{2}(s\mathrm I_H-A_2)^{-1}B_2$. Then, the function $\mathrm I_{\mathbb R^m\times \mathbb R^{p}}-\mathbf H_1(s)\mathbf H_2(s)$ has a proper inverse on some right half-plane.  
\end{itemize}
\end{assumption}
These assumptions lead to the proof that \eqref{eq:spm} is a regular linear system, as explained in \cite{weiss1997dynamic}. 

Due to the extension of $C_1$ (and because $D_1$ appears), the system given in \eqref{eq:spm} is changed as follows:

\begin{equation}
\label{eq:spm}
    \left\{
\begin{aligned}
     &\varepsilon \frac{\diff}{\diff t}z = A_1 z + B_1C_2 w,\\
     & \frac{\diff}{\diff t} w = A_2 w + B_2 (C_{\Lambda,1} z + D_1 C_2 w),\\
     & y = \begin{bmatrix}
         C_{\Lambda,1} z & C_2 w
     \end{bmatrix} \begin{bmatrix}
         z \\ w
     \end{bmatrix} + \begin{bmatrix}
         D_1 & 0
     \end{bmatrix} \begin{bmatrix}
         C_2 w \\ C_1 z
     \end{bmatrix},\\
     & z(0) = z_0,\: w(0) = w_0.
\end{aligned}
    \right.
\end{equation}

As explained in the introduction, one of the goals of the singular perturbation technique consists in finding a way to decouple \eqref{eq:spm} into two subsystems: the first (called the reduced order system) would describes the slower dynamics, while the second one (called the boundary layer system) represents the faster dynamics. The singular perturbation consists in supposing that that, as soon as these approximated systems are exponentially stable, then there exists $\varepsilon^*$ such that for all $\varepsilon\in (0,\varepsilon^*)$, the full-system is exponentially stable. Let us compute them (at least formally).

\paragraph{Reduced order system.} Suppose that $\varepsilon=0$. Then, one has 

$$
\bar z = -A_1^{-1} B_1 C_2 w\in H
$$
Plugging this quantity in the dynamics of $w$ instead of $z$, one obtains:

\begin{equation}
\label{eq:ros}
\left\{
\begin{aligned}
&\frac{\diff}{\diff t} w = (A_2-B_2(C_{\Lambda,1} A_1^{-1}B_1C_2 - D_1 C_2))w,\\
& w(0) = w_0,
\end{aligned}
\right.
\end{equation}
which corresponds to the reduced order system, where $C_{\Lambda,1}$ is the extension of $C_1$ and $D_1$ is given by:
$$
\lim_{\lambda \rightarrow + \infty} \mathbf H_1(\lambda)=D_1,
$$
with $\lambda$ a real value. 

One can prove easily that the operator $B_2(C_{\Lambda,1} A_1^{-1}B_1 - D_1)C_2\in \mathbb R^{n\times n}$. Therefore, it is obvious that there exists a unique solution to \eqref{eq:ros}. We furthermore suppose the following property.

\begin{assumption}
\label{ass:hurwitz-spm}
The matrix $\tilde A_2 := A_2-B_2(C_{\Lambda,1} A_1^{-1}B_1 - D_1)C_2$ is Hurwitz.
\end{assumption}
It means in particular that there exists a skew-symmetric matrix $P_2\in\mathbb R^{n\times n}$ and a positive number $\mu$ such that 
\begin{equation}
\label{eq:Lyap-hurwitz}
2\langle P_2\tilde A_2 w,w\rangle_{\mathbb R^{n\times n}} \leq - \mu \Vert w\Vert^2_{\mathbb R^n} 
\end{equation}

\paragraph{Boundary layer system.} Consider $\tau =\frac{t}{\varepsilon}$. Therefore, the dynamics of $$\tilde z = z + A_1^{-1} B_1 C_2 w$$ satisfies the following:

\begin{equation}
\frac{\diff }{\diff \tau} \tilde z = A_1(z + A_1^{-1}B_1 C_2 w) + \varepsilon \frac{\diff}{\diff t} A_1^{-1} B_1 C_2 w.
\end{equation}
Taking $\varepsilon=0$, one obtains that the boundary layer system can be written as follows:

\begin{equation}
\label{eq:ros}
\left\{
\begin{aligned}
&\frac{\diff}{\diff t} \tilde z = A_1 \tilde z,\\
& \tilde z(0) = \tilde z_0.
\end{aligned}
\right.
\end{equation}

We further consider the "disturbed" version of \eqref{eq:ros} with the output $y=C_1z$. Therefore, one has:
\begin{equation}
\label{eq:ros-disturbed}
\left\{
\begin{aligned}
&\frac{\diff}{\diff t} \tilde z = A_1 \tilde z + \tilde B_1 u,\\
& y_1(t) = C_1 z,\\
& \tilde z(0) = \tilde z_0,
\end{aligned}
\right.
\end{equation}
where $\tilde B_1$ is an admissible control operator for $\mathbb T$ (for instance, $B_1$).
\begin{assumption}
\label{ass:coercive-spm}
Assume that there exists a coercive-ISS Lyapunov functional $V:H\rightarrow \mathbb R_+$ for \eqref{eq:ros-disturbed} given by $\langle Pz,z\rangle_H$ with $P\in \mathcal L(H)$. This Lyapunov functional satisfies, for all $z\in D(A_1)$:
$$
\frac{\diff}{\diff t} V(z) \leq - \alpha_1 \Vert z\Vert_H^2 + \alpha_2 \Vert u\Vert^2_{\mathbb R^p} - \alpha_3 \Vert y_1\Vert^2_{\mathbb R^m},
$$
where $\alpha_1,\alpha_2,\alpha_3>0$. 
\end{assumption}

\begin{remark}
Note that, according to Section \ref{sec:coercive}, once one has an ISS-Lyapunov functional for the system \eqref{eq:ros-disturbed} (without output), then, as soon as the pair $(A_1,C_1)$ is exactly observable, the ISS Lyapunov functional can be made coercive. The converse is also true, meaning that, if we suppose that Assumption \ref{ass:coercive-spm} is satisfied, then the pair $(A_1,C_1)$ is automatically exactly observable. 
\end{remark}

\paragraph{Singular perturbation result.} We are now in position to state our result.

\begin{theorem}
\label{thm:spm}
Suppose that Assumptions \ref{ass:coercive-spm} and \ref{ass:hurwitz-spm} are satisfied. Therefore, there exists $\varepsilon^*$ such that, for every $\varepsilon\in (0,\varepsilon^*)$, the full-system \eqref{eq:spm} is globally exponentially stable.
\end{theorem}

The proof is based on the introduction of a suitable change of variable, namely:
$$
\tilde z = z + A_1^{-1} B_1 C_2 w,
$$
which satisfies the following differential equations:
\begin{equation}
\label{eq:spm-coordinates}
\left\{
\begin{aligned}
&\varepsilon\frac{\diff}{\diff t} \tilde z(t) = A_1 \tilde z(t) + \varepsilon A_1^{-1}B_1 (C_2 \tilde A_2 w(t) + C_2 B_2 C_{\Lambda,1}\tilde z(t))\\
& \frac{\diff}{\diff t} w(t) = \tilde A_2 w(t) + B_2 C_{\Lambda,1} \tilde z(t),\\
& \tilde z(0) = z_0,\: w(0) = w_0,
\end{aligned}
\right.
\end{equation}
where we recall that $\tilde A_2$ is defined as $\tilde A_2 := A_2-B_2(C_{\Lambda,1} A_1^{-1}B_1 - D_1)C_2$.

Now, we are ready to provide the proof of Theorem \ref{thm:spm}:

\begin{proof}
We consider the initial conditions $(z_0,w_0)\in D(A_1)\times \mathbb R^n$ so that the trajectories of \eqref{eq:spm-coordinates} have a suitable regularity (namely, the solutions are strong) to consider time derivative of the Lyapunov functionals that will be introduced later. We consider the Lyapunov functional:
$$
W(z,w):= \varepsilon V(z) + \langle P_2(w-\varepsilon \mathcal Mz),w-\varepsilon \mathcal Mz\rangle_{\mathbb R^n}, 
$$
with $V$ defined in Assumption \eqref{ass:coercive-spm}, $P_2$ defined in \eqref{eq:Lyap-hurwitz}, and where $M$ is given by:
\begin{equation}
\label{eq:M}
    \mathcal M:=B_2C_1A_1^{-1}\in\mathcal L(H,\mathbb R^n).
\end{equation}
One can prove that the function $W$ is coercive, as stated in Lemma \ref{lem:bounds-Lyap} given in the Appendix. Therefore, using the property given in \eqref{ass:coercive-spm} with $y=C_1 \tilde z$, $\tilde B_1 = A^{-1}B_1$ as the admissible operator, and $u(t)=C_2 \tilde A_2 w(t) + C_2 B_2 C_{\Lambda,1}\tilde z(t)$, one obtains that, for all $z\in D(A_1)$ and for all $w\in \mathbb R^n$

\begin{equation}
\begin{aligned}
    \frac{\diff }{\diff t} V(\tilde z) \leq &- \alpha_1 \Vert \tilde z\Vert^2_H + \alpha_2 \varepsilon (2\Vert C_2 \tilde A_2\Vert^2_{\mathcal L(\mathbb {R}^n,\mathbb{R}^p)} \Vert w\Vert^2_{\mathbb R^n} +2 \Vert C_2 B_2\Vert^2_{\mathcal{L}(\mathbb R^m,\mathbb R^p)} \Vert C_{1}\tilde z\Vert_{\mathbb R^p}^2)\\
    &- \alpha_3 \Vert C_1 \tilde z\Vert^2_{\mathbb R^p}
    \end{aligned}
\end{equation}
Moreover, for all $z\in D(A_1)$ and for all $w\in \mathbb R^n$:
\begin{equation}
\begin{aligned}
\frac{\diff}{\diff t} \langle P_2(w-\varepsilon\mathcal M\tilde z),w-\varepsilon&\mathcal M z\rangle_{\mathbb R^n} = 2 \langle P_2(\tilde A_2 w + B_2 C_1 \tilde z,w-\varepsilon \mathcal M z\rangle_{\mathbb R^n}\\
&  + 2 \langle P_2\mathcal M(A_1\tilde  z + \varepsilon A_1^{-1}B_1 (C_2 \tilde A_2 w + C_2 B_2 C_{\Lambda,1}\tilde z),w-\mathcal M \tilde z\rangle_{\mathbb R^n}
\end{aligned}
\end{equation}
Recalling that the definition of $\mathcal M$ given in \eqref{eq:M}, and using the notation $$U(\tilde z,w):= \langle P_2(w-\varepsilon\mathcal M\tilde z),w-\varepsilon\mathcal M z\rangle_{\mathbb R^n},$$ one therefore obtains, for all $(z,w)\in D(A_1)\times\mathbb R^n$:
\begin{equation}
\begin{aligned}
\frac{\diff}{\diff t} U(\tilde z,w) = &2 \langle P_2(\tilde A_2 w,w-\varepsilon \mathcal M \tilde z\rangle_{\mathbb R^n}\\
 & + 2 \varepsilon\langle P_2\mathcal M A_1^{-1}B_1 (C_2 \tilde A_2 w + C_2 B_2 C_{\Lambda,1}\tilde z),w-\mathcal M\tilde  z\rangle_{\mathbb R^n}\\
 \leq & -\mu \Vert w \Vert^2_{\mathbb R^n} + \varepsilon \Vert P_2 \tilde A_2\Vert^2_{\mathcal{L}(\mathbb R^n)} \Vert w\Vert^2_{\mathbb R^n} \\
&+ \varepsilon \Vert P_2\Vert^2_{\mathcal L(\mathbb R^n)} \Vert \mathcal M\Vert^2_{\mathcal L(H,\mathbb R^n)} \Vert \tilde z\Vert_H^2\\
& + 2 \varepsilon \Vert  P_2 \mathcal M A_1^{-1}B_1 C_2 \tilde A_2\Vert_{\mathcal L(\mathbb R^n)} \Vert w\Vert^2_{\mathcal L(\mathbb R^n)}\\
&  + \varepsilon \Vert P_2 \mathcal M A_1^{-1}B_1 C_2 B_2\Vert^2_{\mathcal L(\mathbb R^n,\mathbb R^p)} \Vert C_1 \tilde z\Vert^2_H\\
& + \varepsilon \Vert  w\Vert_{\mathbb R^n}^2 + \varepsilon \Vert  P_2 \mathcal M A_1^{-1}B_1 C_2\tilde A_2\Vert_{\mathcal L(\mathbb R^n)}^2 \Vert w\Vert^2_{\mathcal L(\mathbb R^n)} \\
& + \varepsilon \Vert \mathcal M\Vert^2_{\mathcal L(H,\mathbb R^n)} \Vert \tilde z\Vert_H^2\\
&+ \varepsilon \Vert P_2 \mathcal M A_1^{-1} B_1 C_2 B_2\Vert^2_{\mathcal L(\mathbb R^p,\mathbb R^n)} \Vert C_1 \tilde z\Vert^2_{\mathbb R^p}\\
& + \varepsilon \Vert \mathcal M\Vert^2_{\mathcal{L}(H,\mathbb R^n)} \Vert \tilde z\Vert_H^2
\end{aligned}
\end{equation}
Using the following notation,
\begin{equation}
\begin{aligned}
    a_1 = &2 \alpha_2 2\Vert C_2 \tilde A_2\Vert^2_{\mathcal L(\mathbb {R}^n,\mathbb{R}^p)} + \Vert P_2 \tilde A_2\Vert^2_{\mathcal{L}(\mathbb R^n)} + 2 \Vert  P_2 \mathcal M A_1^{-1}B_1 C_2 \tilde A_2\Vert_{\mathcal L(\mathbb R^n)} \\
    &+ 1 + \Vert  P_2 \mathcal M A_1^{-1}B_1 C_2 \tilde A_2\Vert_{\mathcal L(\mathbb R^n)}^2,\\
    a_2 = & 2\left(\alpha_2 \Vert C_2B_2\Vert^2_{\mathcal L(\mathbb R^m,\mathbb R^p)} +   \Vert P_2 \mathcal M A_1^{-1}B_1 C_2 B_2\Vert^2_{\mathcal L(\mathbb R^n,\mathbb R^p)}\right),\\
    a_3  = & \Vert P_2\Vert^2_{\mathcal L(\mathbb R^n)} \Vert \mathcal M\Vert^2_{\mathcal L(H,\mathbb R^n)} + 2 \Vert \mathcal M\Vert^2_{\mathcal L(H,\mathbb R^n)},
\end{aligned}
\end{equation}
one therefore has, for any $(z,w)\in D(A_1)\times \mathbb R^n$:
\begin{equation}
\frac{\diff }{\diff t} W(z,w) \leq -(\mu - \varepsilon a_1) \Vert w\Vert^2_{\mathbb R^n} - (\alpha_3 -\varepsilon a_2) \Vert C_1 z\Vert^2_{\mathbb R^p} - (\alpha_1 - \varepsilon a_3)\Vert z\Vert_H^2.
\end{equation}
Therefore, for any $\varepsilon$ such that
$$
\varepsilon <\varepsilon^*:= \min\left(\frac{\mu}{a_1},\frac{\alpha_3}{a_2}, \frac{\alpha_1}{a_3}\right),
$$
one obtains that
\begin{equation}
    \frac{\diff}{\diff t} W(\tilde z,w) \leq - \tilde \mu \Vert w\Vert_{\mathbb R^n} - \tilde \alpha_3 \Vert C_1 \tilde z\Vert^2_{\mathbb R^p} - \tilde \alpha_1 \Vert \tilde z\Vert_H^2\leq -  \tilde \mu \Vert w\Vert_{\mathbb R^n} - \tilde \alpha_1 \Vert \tilde z\Vert_H^2.
\end{equation}
Using the coercivity of $W$ (proved in Lemma \ref{lem:bounds-Lyap}), one can can prove the desired result for any strong solutions. The same result can be deduced by a density argument for weak solutions.
\end{proof}

\begin{remark}
This result might be seen as weak, since it does tackle the case of coupled infinite-dimensional systems. However, as it is well illustrated in \cite{tang2017stability,cerpa-prieur2017}, there exist counter-examples to the singular perturbation method for slow infinite-dimensional systems coupled with fast ODEs. Therefore, we cannot hope having a general result concerning coupled infinite-dimensional systems, unless some strong assumptions are given. 
\end{remark}

\subsection{Output regulation: the PI controller case}

Given the real Hilbert space $H$, let us consider now the system:
\begin{equation}
\label{eq:PI}
\left\{
\begin{aligned}
& \frac{\diff}{\diff t} z = A z + Bw,\\
& \frac{\diff}{\diff t} w = k(y-r),\\
& y = Cz\\
& z(0)=z_0,\: w(0) = w_0.
\end{aligned}
\right.
\end{equation}
where $A:D(A) \subset H\rightarrow H$ is the generator of a strongly continuous semigroup $\mathbb T$, $B\in\mathcal L(\mathbb R,H_{-1})$, $C\in \mathcal L(H_1,\mathbb R)$, $(A,B,C)$ is supposed to be a regular triple, and $w,r\in \mathbb R$\footnote{The result is also true when consider vectors instead of scalars. We chose this framework in order to make easier the reading of the paper.}. This system is an example of an infinite-dimensional system controlled with a PI controller (where $w$ denotes the state of the integrator and $r$ the reference). As proved in \cite{paunonen2014internal}, the system \eqref{eq:PI} is well-posed and admits a unique generalized solution. For this, we need to consider the extension of $C$, denoted by $C_\Lambda$. Hence, one has
\begin{equation}
    \label{eq:PI-extended}
 \left\{
\begin{aligned}
& \frac{\diff}{\diff t} z = A z + Bw,\\
& \frac{\diff}{\diff t} w = k(y-r),\\
& y = C_\Lambda z + Dw\\
& z(0)=z_0,\: w(0) = w_0.
\end{aligned}
\right.   
\end{equation}

To put \eqref{eq:PI} in the same form than \eqref{eq:spm}, we compute its equilibrium points, denoted by $z_e$ and $w_e$:
\begin{equation}
C_\Lambda z_e + Dw_e  = r,\: z_e = - A^{-1} Bw_e
\end{equation}
Therefore, one has:
\begin{equation}
y_e = r,\: w_e= -(C_\Lambda A^{-1} B -D)r. 
\end{equation}
In order to make sure that $w^*$ exists, one has to suppose the following:
$$
(C_\Lambda A^{-1} B -D)\neq 0,
$$
which is a classical assumption in output regulation. It says that, at $s=0$, the transfer function $\mathbf H(s) = C_\Lambda (\mathrm sI_H-A)^{-1}B + D$ is different from $0$.

Consider now the following change of coordinates:
$$
\tilde z = z-z_e,\: \tilde w= w-w_e.
$$
These variables satisfy the dynamics given as follows:
\begin{equation}
    \label{eq:PI-extended-coordinates}
 \left\{
\begin{aligned}
& \frac{\diff}{\diff t} \tilde z = A \tilde z + B\tilde w,\\
& \frac{\diff}{\diff t} \tilde w = k \tilde y,\\
& \tilde y = y - r\\
& z(0)=z_0,\: w(0) = w_0.
\end{aligned}
\right.   
\end{equation}
Note moreover that, as noted in \cite{lorenzetti2022saturating}, the PI controller is referred to as a low-gain technique (in contrast with high-gain techniques such as the backstepping \cite{krstic_smyshlyaev_backstepping}). This implies that $k$ can be seen as a small scalar. Therefore, performing the change of time-scale $\tau:=k t$, one obtains:
\begin{equation}
    \label{eq:PI-extended-spm}
 \left\{
\begin{aligned}
& k\frac{\diff}{\diff \tau} \tilde z = A \tilde z + B\tilde w,\\
& \frac{\diff}{\diff \tau} \tilde w = \tilde y,\\
& \tilde y = y - r,\\
& z(0)=z_0,\: w(0) = w_0.
\end{aligned}
\right.   
\end{equation}
which means that the PI controller problem given in \eqref{eq:PI} can be seen as a singular perturbation problem as the one given in \eqref{eq:spm}. 

\paragraph{Reduced order system.} The reduced order system can be obtained by taking $k=0$. Then, one has:
$$
\bar z = - A^{-1} B\tilde w,
$$
which implies that
\begin{equation}
    \label{eq:ros-PI}
\frac{\diff}{\diff \tau } \tilde w = (-C_\Lambda A^{-1} B + D)\tilde w,
\end{equation}
meaning that we need to assume the following for the reduced order system:
\begin{assumption}
\label{ass:ros-PI}
The scalar $-C_\Lambda A^{-1} B + D$ is negative.
\end{assumption}

\paragraph{Boundary layer system.} Now, consider the initial time $t$ and the $z_e=\tilde z + A^{-1} B \tilde w$. Then, we have 
\begin{equation}
\label{eq:bls-PIk}
\frac{\diff}{\diff t} \tilde z  = A\tilde z + k(A^{-1}B C_\Lambda - D)\tilde w
\end{equation}
Taking $k=0$, one obtains that
\begin{equation}
\label{eq:bls-PI}
\left\{
\begin{aligned}
&\frac{\diff}{\diff t} \tilde z = A\tilde z,\\
& y= Cz,\\
&\tilde z(0) = z_0,
\end{aligned}
\right.
\end{equation}

As in the singular perturbation case, let us consider an admissible operator $\tilde B$ and an input $u\in L^2_{\mathrm loc}([0,\infty))$.
\begin{equation}
\label{eq:bls-PI-u}
\left\{
\begin{aligned}
&\frac{\diff}{\diff t} \tilde z = A\tilde z + \tilde B u,\\
& y= Cz,\\
&\tilde z(0) = z_0,
\end{aligned}
\right.
\end{equation}
Therefore, we assume the following:
\begin{assumption}
\label{ass:coerciveISSPI}
We assume that the pair $(A,B,C)$ is regular and that there exists a coercive ISS Lyapunov function $V:=\langle Pz,z\rangle_H$ satisfying along the strong solution to \eqref{eq:bls-PI-u}:
\begin{equation}
    \frac{\diff}{\diff t} V(z) \leq - \alpha_1 \Vert z\Vert_H^2 + \alpha_2 |u|^2 - \alpha_3 |y|^2,
\end{equation}
with $\alpha_1,\alpha_2,\alpha_3$ positive numbers. 
\end{assumption}

Under these assumptions, one can easily apply Theorem \ref{thm:spm} to find bounds on $k$ in order to stabilize (and output regulate) \eqref{eq:PI}.

\begin{example}[Korteweg-de Vries equation (continued)]
Consider the KdV equation given in Example \ref{example:kdv}. Let us see how one can regulate it. We therefore add in KdV equation presented in Example \ref{example:kdv} an integrator, whose state is denoted by $w$.
\begin{equation}
    \left\{
\begin{aligned}
&z_t + z_x + z_{xxx} = 0,\: &\text{ on } &(t,x)\in [0,\infty)\times [0,L]\\
& z(t,0) = z(t,L) = 0,\: &\text{ on } &t\in [0,\infty),\\
&z_x(t,L) = w(t),\: &\text{ on } & t\in [0,\infty),\\
& y(t) = z_x(t,0),\: &\text{ on } & t\in [0,\infty),\\
& \frac{\diff}{\diff t} w (t)=k( y(t) - r),\: & \text{ on } &t\in [0,\infty),\\
& z(0,x)= z_0(x),\: w(0) = w_0,\: & \text{ on } &x\in [0,L].
\end{aligned}
    \right.
\end{equation}
The equilibrium points are given by:
\begin{equation}
\label{eq:bvp-KdV}
\left\{
\begin{aligned}
&z^\prime_e(x) + z^{\prime \prime\prime}_e(x) =0,\\
&z_e(0) = z_e(L) = 0,\\
& z_e^\prime(L) = w_e,
\end{aligned}
\right.
\end{equation}
and
\begin{equation}
\label{eq:equilibrium}
z_e^\prime(0) = r.
\end{equation}
One can find an explicit solution to \eqref{eq:bvp-KdV} together with \eqref{eq:equilibrium} given by
$$
z_e(x) = 2r\frac{\sin(\frac{x}{2}) \sin(\frac{L-x}{2})}{\sin(\frac{L}{2})}.
$$
Therefore, one can deduce $w_e$ from \eqref{eq:bvp-KdV}. After the same change of coordinates than in the section presented before, one obtains:
\begin{equation}
    \left\{
\begin{aligned}
&\tilde z_t + \tilde z_x + \tilde z_{xxx} = 0,\: &\text{ on } &(t,x)\in [0,\infty)\times [0,L]\\
& \tilde z(t,0) = \tilde z(t,L) = 0,\: &\text{ on } &t\in [0,\infty),\\
&\tilde z_x(t,L) = \tilde w(t),\: &\text{ on } & t\in [0,\infty),\\
& \tilde y(t) = \tilde z_x(t,0),\: &\text{ on } & t\in [0,\infty),\\
& \frac{\diff}{\diff t} w (t)=k\tilde y(t),\: & \text{ on } &t\in [0,\infty),\\
& z(0,x)= z_0(x),\: w(0) = w_0,\: & \text{ on } &x\in [0,L].
\end{aligned}
    \right.
\end{equation}
We change the time scale with $\tau = k t$. Therefore, one has:
\begin{equation}
\label{eq:spmKdV}
    \left\{
\begin{aligned}
&k \tilde z_\tau + \tilde z_x + \tilde z_{xxx} = 0,\: &\text{ on } &(t,x)\in [0,\infty)\times [0,L]\\
& \tilde z(\tau,0) = \tilde z(\tau,L) = 0,\: &\text{ on } &\tau\in [0,\infty),\\
&\tilde z_x(\tau,L) = \tilde w(\tau),\: &\text{ on } & \tau\in [0,\infty),\\
& \tilde y(\tau) = \tilde z_x(\tau,0),\: &\text{ on } & \tau\in [0,\infty),\\
& \frac{\diff}{\diff \tau } w (\tau)=\tilde y(\tau),\: & \text{ on } &\tau\in [0,\infty),\\
& z(0,x)= z_0(x),\: w(0) = w_0,\: & \text{ on } &x\in [0,L].
\end{aligned}
    \right.
\end{equation}
We are therefore in the situation described by \eqref{eq:spm}. Now, let us check whether the system is regular. We have proved in Example \ref{example:kdv} that the subsystem given by the KdV equation is regular (i.e., the system with the unknown $z$). The operators $A,B,C$ are defined in Example \ref{example:kdv} and recall that the associated transfer function $\mathbf H_1(s)=1$. The associated transfer function associated to the system given by the integrator is $\mathbf H_2(s) = \frac{1}{s}$. It is clear that $1 - \mathbf H_1(s) \mathbf H_2(s)$ admits a proper inverse on some right half plane, given by:
$$
(1 - \mathbf H_1(s) \mathbf H_2(s))^{-1} = \frac{s}{s-1},
$$
which shows that the system represented by \eqref{eq:spmKdV} is regular, i.e. there exists a unique strong (resp., generalized) solution to \eqref{eq:spmKdV} for any initial condition $(z_0,w_0)\in D(S)\times \mathbb R$ (resp., for any initial condition $(z_0,w_0)\in H\times \mathbb R$), where $D(S)$ is given by
$$
D(S) := \lbrace (z,w)\in H^3(0,L)\times \mathbb R\mid z(0) = z(L) = 0,\: z^\prime(L) = w\rbrace. 
$$

We can now compute the subsystems, namely the reduced order system and the boundary layer system. 

\paragraph{Reduced order system.} Take $k=0$ and observer that, for all $\tau\geq 0$ (which is seen as a fixed parameter):
\begin{equation}
    \left\{
\begin{aligned}
&\tilde z_x + \tilde z_{xxx} = 0,\: \text{ on } x\in [0,L]\\
&\tilde z(\tau,0) = \tilde z(\tau,L) = 0,\\
& \tilde z(\tau, L) = \tilde w(\tau),
\end{aligned}
    \right.
\end{equation}
where it can be shown that:
$$
\tilde z(\tau,x) = - 2\tilde w(\tau)\frac{\sin(\frac{x}{2})\sin(\frac{L-x}{2})}{\sin(\frac{L}{2})}.
$$
Therefore, the reduced order system is given by:
\begin{equation}
\left\{
\begin{aligned}
& \frac{\diff}{\diff \tau} \tilde w(\tau) = - \tilde w(\tau),\\
& \tilde w(0) = w_0,
\end{aligned}
\right.
\end{equation}
which is clearly stable. Assumption \ref{ass:ros-PI} is therefore satisfied.

\paragraph{Boundary Layer system.} Consider now, for all $(\tau,x)\in [0,\infty)\times [0,L]$, $$\bar z(\tau,x) = \tilde z(\tau,x) + 2\tilde w(\tau)\frac{\sin(\frac{x}{2})\sin(\frac{L-x}{2})}{\sin(\frac{L}{2})},$$ which is the difference between the state $\tilde z$ and the equilibrium state expressed with $\tilde w$. Therefore, taking $t=\frac{\tau}{k}$, and then $k=0$, one obtains that
\begin{equation}
\left\{
\begin{aligned}
& \bar z_t + \bar z_x + \bar z_{xxx} =0, &\text{ on } &[0,\infty)\times [0,L],\\
& \bar z(t,0) = \bar z(t,L) = 0, &\text{ on } &[0,\infty),\\
& \bar z_x(t,L) = 0, & \text{ on } &[0,\infty) \times [0,L],\\
& \bar z(0,x) = \tilde z_0(x) + 2\tilde w(0)\frac{\sin(\frac{x}{2})\sin(\frac{L-x}{2})}{\sin(\frac{L}{2})}, & \text{ on } & [0,L],
\end{aligned}
\right.
\end{equation}
which is exponentially stable as soon as $L\notin \mathcal N$ (see e.g., \cite{cerpa2013control,rosier1997kdv}), where we recall that $\mathcal N$ is defined in Example \ref{example:kdv}. Moreover, it can be shown that the pair $(A,C)$, expressed in Example \ref{example:kdv}, is exactly observable as soon as $L\notin \mathcal N$ (which we will suppose in the following). We even know a coercive ISS-Lyapunov functional (called $V$), obtained in \cite{balogoun2021iss} and also used in \cite{marx2024singular}. Assumption \ref{ass:coerciveISSPI} is therefore satisfied. 

Therefore, Theorem \ref{thm:spm} can be applied, and we can therefore find a bound on $k$ so that the system \ref{eq:spmKdV} is exponentially stable.
\end{example}

\begin{remark}
Building a PI controller for the KdV equation was the goal of \cite{balogoun2021iss}, and the bounds obtained in the latter paper are surely more precise than the ones expressed in this paper. We just wanted to emphasize on the link between singular perturbation and output regulation, and how useful a coercive ISS Lyapunov functional can be for such problems.
\end{remark}

\section{Conclusion}

\label{sec:conclusion}

This paper has proposed a method which allows to make an ISS non-coercive Lyapunov functional coercive, under the condition that there exists an output operator $C$ such that the pair $(A,C)$ is exactly observable. The technique is related to the strictification technique presented in \cite{praly2019observers} and is inspired mainly by \cite{lorenzetti2022saturating}. However, adding the input needed us to consider a subclass well-posed, namely the regular linear systems, that are very useful when one wants to extend known finite-dimensional techniques to infinite-dimensional systems.

We believe that the construction of such coercive Lyapunov functional might be useful for other problems than the ones proposed in the papers, such as the construction of observers. We also feel that this functional could be useful for linear systems admitting some nonlinear dampings, such as the ones presented in \cite{marx2021forwarding,chitour2019p,chitour2020one,singh2023second,marx2025impedance}.

\section*{Appendix}
\label{sec:appendix}

\subsection*{Coercivity of $W$}
This short appendix is devoted to the proof of a result used in the proof of Theorem \ref{thm:spm}, which uses a Lyapunov functional. Here is the statement of such a result. 
\begin{lemma}
\label{lem:bounds-Lyap}
Suppose that Assumptions \ref{ass:coercive-spm} and \ref{ass:hurwitz-spm} are satisfied. Consider the Lyapunov functional:
$$
W(z,w):= \varepsilon V(z) + \langle P_2(\mathcal M z-w),\mathcal Mz - w\rangle_{\mathbb R^n}, 
$$
where $P_2$ is given by \eqref{eq:Lyap-hurwitz} and $\mathcal{M}$ given by \eqref{eq:M}. Then, there exists $\bar \nu,\underline \nu>0$ such that
\begin{equation}
  \underline \nu \left(\Vert z\Vert^2_{H} + \Vert w \Vert_{\mathbb R^n}^2\right)\leq  W(z,w) \leq \overline\nu \left(\Vert z\Vert^2_{H} + \Vert w \Vert_{\mathbb R^n}^2\right),
\end{equation}
where $$\overline \nu := \max (\varepsilon \overline p+ \varepsilon^2 \Vert \mathcal M\Vert_{\mathcal L(H,\mathbb R^n)}^2 \Vert P_2\Vert_{\mathcal L(W)},\Vert P_2\Vert_{\mathcal L(\mathbb R^n)}),$$ and $$\underline \nu:= \min \left(\frac{\underline p\varepsilon}{2},\frac{1}{2}\frac{\underline p\varepsilon}{\varepsilon^2 \beta \Vert \mathcal M\Vert^2_{\mathcal L(H,\mathbb R^n)} + a_1 \varepsilon}\right),$$ with $\beta$ the eigenvalue of lowest value of $P_2$.
\end{lemma}

\begin{proof}
Using Young's Lemma and the fact that Assumption \ref{ass:coercive-spm} is satisfied, it is clear that
\begin{equation}
\begin{split}
W(z,w) \leq &\Vert P_2\Vert_{\mathcal L(\mathbb R^n)}\Vert w\Vert^2_{\mathbb R^n} + \varepsilon \overline p \Vert z\Vert_H^2\\
&+ \varepsilon^2 \Vert \mathcal M\Vert_{\mathcal L(H,\mathbb R^n)}^2 \Vert P_2\Vert_{\mathcal L(\mathbb R^n)} \Vert z\Vert_Z^2.
\end{split}
\end{equation}
Moreover, using again Young's Lemma and the fact that Assumptions \ref{ass:coercive-spm} and \ref{ass:hurwitz-spm} are satisfied, one obtains:
\begin{equation*}
\begin{split}
    W(z,w) \geq & \underline p \varepsilon \Vert z\Vert^2_H \\
    &+ \beta(\varepsilon^2\Vert\mathcal M z\Vert_{\mathbb R^n}^2 + \Vert w\Vert_{\mathbb R^n}^2 - 2\langle \varepsilon \mathcal Mz,w\rangle_{\mathbb R^n})\\
\geq & \underline p \varepsilon \Vert z\Vert^2_H + \frac{\beta \varepsilon^2}{2}\left(1-\frac{1}{\delta}\right)\Vert\mathcal Mz\Vert_{\mathbb R^n}^2 \\
&+ \frac{\beta }{2}(1-\delta)\Vert w\Vert_{\mathbb R^n}^2
\end{split}
\end{equation*}
Select $\delta:=\frac{\beta \varepsilon^2\Vert \mathcal M\Vert^2_{\mathcal L(H,\mathbb R^n)}}{\beta \varepsilon^2\Vert M\Vert^2_{\mathcal L(H,\mathbb R^n)} + \overline p \varepsilon}$. Then, one has $1-\frac{1}{\delta}<0$. Moreover,
\begin{equation*}
\begin{split}
    W(z,w) \geq & \overline p \varepsilon \Vert z\Vert_Z^2 \\
    &- \frac{1}{2}\left(\frac{\overline p \varepsilon}{\beta \varepsilon^2 \Vert \mathcal M\Vert^2_{\mathcal{L}(H,\mathbb R^n)}}\right)\varepsilon^2 \Vert M\Vert^2_{\mathcal{L}(H,\mathbb R^n)} \Vert z\Vert^2_Z \\
    &+ \frac{1}{2}\frac{\overline p\varepsilon}{\varepsilon^2 \beta \Vert \mathcal M\Vert^2_{\mathcal{L}(H,\mathbb R^n)}+ \overline p\varepsilon} \Vert w\Vert_{\mathbb R^n}^2,
\end{split}
\end{equation*}
meaning in particular that the statement of Lemma \ref{lem:bounds-Lyap} holds true. This concludes the proof.
\end{proof}

%\begin{acknowledgements}
%If you'd like to thank anyone, place your comments here
%and remove the percent signs.
%\end{acknowledgements}

% Authors must disclose all relationships or interests that 
% could have direct or potential influence or impart bias on 
% the work: 
%
% \section*{Conflict of interest}
%
% The authors declare that they have no conflict of interest.

% BibTeX users please use one of
%\bibliographystyle{spbasic}      % basic style, author-year citations
\bibliographystyle{spmpsci}      % mathematics and physical sciences
\bibliography{bibsm}   % name your BibTeX data base

\end{document}